\documentclass[11pt]{amsart}
\usepackage{amsmath,graphicx}

\newtheorem{thm}{Theorem}[section]

\newtheorem{lem}[thm]{Lemma}%[section]
\newtheorem{cor}[thm]{Corollary}%[section]
\theoremstyle{remark}
\newtheorem{remark}[thm]{Remark}%[section]
\theoremstyle{definition}

\theoremstyle{plain}

\newcommand{\Z}{{\mathbb{Z}}}
\newcommand{\Zf}{\Z^2_{\mathrm{vis}}}

\newcommand{\R}{{\mathbb{R}}}

\newcommand{\N}{{\mathbb{N}}}

\newcommand{\TT}{{\mathcal{T}}}
\newcommand{\FF}{{\mathcal{F}}}

\newcommand{\ar}{\operatorname{area}}
\newcommand{\lh}{\operatorname{length}}

\numberwithin{equation}{section}

\begin{document}

\title[]{On the index of Farey sequences}

\author[]{Florin P. Boca, Radu N. Gologan and Alexandru Zaharescu}

\address{FPB and AZ: Department of Mathematics, University of
Illinois at Urbana-Champaign, Urbana, IL 61801, USA;
fboca@math.uiuc.edu; zaharesc@math.uiuc.edu}

\address{Institute of Mathematics of the Romanian Academy,
P.O. Box 1-764, RO-014700 Bucharest, Romania}

\address{RNG: Institute of Mathematics of the Romanian Academy,
P.O. Box 1-764, RO-014700 Bucharest, Romania;
Radu.Gologan@imar.ro}

\thanks{Research partially supported by ANSTI grant C6189/2000}

%\date{March 12, 2002}

\begin{abstract}
We prove some asymptotic formulae concerning the distribution of
the index of Farey fractions of order $Q$ as $Q\rightarrow \infty$.
\end{abstract}

\maketitle

\section{Introduction}

Let $\FF_Q =\{ \gamma_1,\dots,\gamma_{N(Q)} \}$ denote the Farey
sequence of order $Q$ with $\gamma_1 =1/Q <\gamma_2 <\cdots <
\gamma_{N(Q)} =1$. This sequence is extended by
$\gamma_{i+N(Q)}=\gamma_i+1$, $1\leq i\leq N(Q)$. It is
well-known that
\begin{equation*}
N(Q)=\sum\limits_{j=1}^Q \varphi (j)
=\frac{3 Q^2}{\pi^2} +O(Q\log Q).
\end{equation*}

\begin{figure}[ht]
\begin{center}
\includegraphics*[scale=0.7, bb=0 0 250 250]{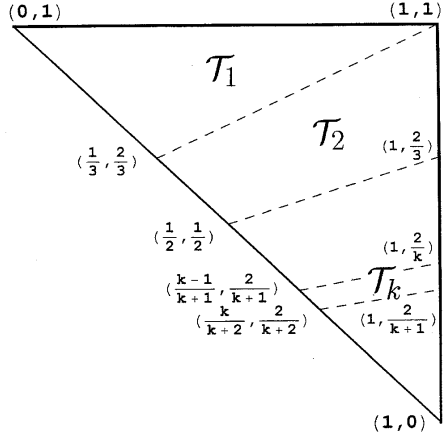}
\end{center}
\caption{The decomposition $\TT =\cup_{k=1}^\infty \TT_k$.}
\label{Figure1}
\end{figure}

For any two consecutive Farey fractions $\gamma_i =a_i /q_i
<\gamma_{i+1}=a_{i+1} / q_{i+1}$,
one has $a_{i+1}q_i-a_i q_{i+1}=1$ and $q_i+q_{i+1}>Q$.
Conversely, if $q$ and $q^\prime$ are two coprime integers
in $\{ 1,\dots,Q\}$ with $q+q^\prime >Q$, then there are
unique $a\in \{ 1,\dots,q\}$ and
$a^\prime \in \{ 1,\dots ,q^\prime \}$ for which
$a^\prime q-aq^\prime =1$, and $a/q <
a^\prime / q^\prime$ are consecutive Farey fractions
of order $Q$. Therefore, the pairs of coprime integers
$(q,q^\prime)$ with $q+q^\prime >Q$ are in one-to-one
correspondence with the pairs of consecutive Farey fractions of
order $Q$. Moreover, the denominator $q_{i+2}$ of $\gamma_{i+2}$
can be easily expressed (cf.\,\cite{HT}) by means of the
denominators of $\gamma_i$ and $\gamma_{i+1}$ as
\begin{equation*}
q_{i+2}=\left[ \frac{Q+q_i}{q_{i+1}} \right] q_{i+1}-q_i .
\end{equation*}
This formula was recently used in the study of various statistical
properties of the Farey fractions \cite{ABCZ}, \cite{BCZ},
\cite{HS}, leading to the definition in \cite{BCZ} of a new and
interesting area-preserving transformation $T$ of the Farey
triangle
\begin{equation*}
\TT =\{ (x,y)\in [0,1]^2 \, ;\, x+y>1\} ,
\end{equation*}
defined by
\begin{equation*}
T(x,y)=\bigg( y,\bigg[ \frac{1+x}{y} \bigg] y-x\bigg)=
\left( y,1-x-y G\bigg( \frac{y}{1+x} \bigg)\right) ,
\end{equation*}
where $G(t)=\{ 1/t \}$ is the Gauss map on the unit interval.
The set $\TT$ decomposes as the union of disjoint sets
(see Figure \ref{Figure1})
\begin{equation*}
\TT_k =\left\{ (x,y)\in \TT \, ;\, \bigg[ \frac{1+x}{y} \bigg]
=k\right\} \qquad k\in \N^* =\{ 1,2,3\dots \} ,
\end{equation*}
and
\begin{equation*}
T(x,y)=(y,ky-x),\qquad (x,y) \in \TT_k .
\end{equation*}

For later use, we also define the upper triangles $\TT_k^\prime$
with vertices at $(1,2/k)$, $\big( (k-1)/(k+1)\big)$ and $\big(
1,2/(k+1)\big)$, $k\geq 2$, and the lower triangles $\TT_k^{\prime
\prime}$ with vertices at $\big( k/(k+2),2/(k+2)\big)$, $\big(
(k-1)/(k+1),2/(k+1)\big)$ and $\big( 1,2/(k+1)\big)$, $k\geq 1$.

For all integer $i\geq 0$, we have
\begin{equation*}
T^i(x,y)=\big( L_i (x,y),L_{i+1}(x,y)\big),
\end{equation*}
with
\begin{equation*}
\begin{split}
& L_{i+1}(x,y)=\kappa_i(x,y)L_i(x,y)-L_{i-1}(x,y) ,\\ &
L_0(x,y)=x,\quad L_1(x,y)=y \\
& \kappa_i(x,y)=\kappa_{i-1} \circ  T(x,y),\quad
\kappa_1(x,y)=[(1+x)/y].
\end{split}
\end{equation*}

It was noticed and extensively used in \cite{BCZ} that
\begin{equation*}
T\bigg( \frac{q}{Q} ,\frac{q^\prime}{Q} \bigg) =
\bigg( \frac{q^\prime}{Q},\frac{q^{\prime \prime}}{Q} \bigg) ,
\end{equation*}
whenever $q$, $q^\prime$ and $q^{\prime \prime}$ are denominators
of three consecutive Farey fractions from $\FF_Q$. This shows
immediately that
\begin{equation}\label{I1}
\kappa_1 \bigg( \frac{q_{i-1}}{Q} ,\frac{q_i}{Q} \bigg) =
\bigg[ \frac{Q+q_{i-1}}{q_i} \bigg]
\end{equation}
and
\begin{equation}\label{I2}
\begin{split}
\kappa_{r+1} \bigg( \frac{q_{i-1}}{Q} ,\frac{q_i}{Q} \bigg) & =
\kappa_1 \circ T^r \bigg( \frac{q_{i-1}}{Q},\frac{q_i}{Q} \bigg) =
\kappa_1 \bigg( \frac{q_{i+r-1}}{Q},\frac{q_{i+r}}{Q} \bigg) \\
&  = \bigg[ \frac{Q+q_{i+r-1}}{q_{i+r}} \bigg] ,\qquad r\in \N .
\end{split}
\end{equation}

We also note that
\begin{equation}\label{I3}
\TT_k =\{ (x,y)\in \TT \, ;\, \kappa_1 (x,y)=k\} .
\end{equation}

In \cite{HS}, the integer
\begin{equation}\label{I4}
\nu_Q (\gamma_i)=\left[ \frac{Q+q_{i-1}}{q_i} \right]
=\kappa_1 \bigg( \frac{q_{i-1}}{Q},\frac{q_i}{Q} \bigg) =
\frac{q_{i+1}+q_{i-1}}{q_i} =\frac{a_{i+1}+a_{i-1}}{a_i}
\end{equation}
was called the index of the fraction $\gamma_i$ in $\FF_Q$,
and various new and interesting results concerning their
distribution were proved, including the striking
closed form formulae
\begin{equation}\label{I5}
\sum\limits_{i} \nu_Q (\gamma_i)=3N(Q)-1
\end{equation}
and
\begin{equation*}
\begin{split}
\sum\limits_{q=1}^Q \# \bigg\{ \gamma_i & =\frac{a_i}{q_i} \, ;\,
q_i=q \ \mbox{\rm and} \ \nu_Q (\gamma_i)=
\bigg[ \frac{2Q+1}{q_i} \bigg] -1\bigg\} \\
& =Q(2Q+1)-N(2Q)-2N(Q)+1 .
\end{split}
\end{equation*}
The following asymptotic formulae were also proved in \cite{HS}:
\begin{equation*}
\sum\limits_{i} \nu_Q (\gamma_i)^2 =
\frac{24}{\pi^2} \, Q^2 \left( \log 2Q -
\frac{\zeta^\prime(2)}{\zeta(2)} -\frac{17}{8}+2\gamma \right)
+O(Q\log^2 Q),
\end{equation*}
\begin{equation}\label{I6}
\begin{split}
& L(Q,k)  =l_k N(Q)+O\bigg( k+\frac{Q\log Q}{k} \bigg) , \\
& U(Q,k) =u_k N(Q)+O\bigg( k+\frac{Q\log Q}{k} \bigg) ,
\end{split}
\end{equation}
where
\begin{equation*}
\begin{split}
& L(Q,k)=\# \bigg\{ \gamma_i \, ;\,
\nu_Q(\gamma_i)=k=\bigg[ \frac{2Q+1}{q_i}\bigg]-1 \bigg\},\\
& U(Q,k)=\# \bigg\{ \gamma_i \, ;\, \nu_Q(\gamma_i)=k
=\bigg[ \frac{2Q+1}{q_i}\bigg] \bigg\} ,
\end{split}
\end{equation*}
and
\begin{equation*}
\begin{split}
& l_k =4\left( \frac{1}{(k+1)^2}-\frac{1}{k+1}+\frac{1}{k+2} \right)
, \\ &  u_k =\begin{cases}
0 & \mbox{\rm if $k=1$},\\
4\Big( \frac{1}{k}-\frac{1}{k+1}-\frac{1}{(k+1)^2} \Big) &
\mbox{\rm if $k\geq 2$.}
\end{cases}
\end{split}
\end{equation*}

In an earlier version of \cite{HS}, it was proved that
\begin{equation}\label{I7}
\sum\limits_{\gamma_i \leq t} \nu_Q (\gamma_i)=3N(Q)t
+O(Q^{3/2+\varepsilon}) \quad
\mbox{\rm whenever $t\in [0,1]$.}
\end{equation}

It was also conjectured that for every $h\in \N^*$,
there is a constant $A(h)$ such that
\begin{equation*}
S_h(Q) =\sum\limits_{i} \nu_Q (\gamma_i) \nu_Q (\gamma_{i+h}) \sim
A(h)N(Q)
\end{equation*}
as $Q\rightarrow \infty$.
In this note, we first prove that this conjecture holds, finding
also the finite constant
\begin{equation*}
A(h)=2\iint\limits_\TT \kappa_1 (s,t)\kappa_{h+1} (s,t) \,
ds\, dt =2\sum\limits_{m,n=1}^\infty mn \ar (T^h\TT_m \cap \TT_n ).
\end{equation*}
All the numbers $A(h)$ are rational and
$A(1)=\frac{192}{35} \approx 5.4857$,
$A(2)=\frac{796727}{90090}
\approx 8.8437$. One can also show that $A(h)\ll 1+\log h$.
It would be interesting to investigate whether $A(h)$ is bounded by
an absolute constant.

We also consider the more general situation where the Farey
fractions belong to a subinterval of $[0,1]$. More precisely, we set
\begin{equation*}
S_{h,t}(Q)=\sum\limits_{\gamma_i \in \FF_Q \cap [0,t]}
\nu_Q (\gamma_i)\nu_Q (\gamma_{i+h}),
\end{equation*}
and prove
\medskip

\begin{thm}\label{ThmA}
{\em (i)} For every integer $h\geq 1$, we have
\begin{equation*}
S_h(Q)= A(h)N(Q) +O_h (Q\log^2 Q)
\end{equation*}
as $Q\rightarrow \infty$.

{\em (ii)} For every integer $h\geq 1$ and every $t\in [0,1]$, we have
for all $\varepsilon >0$,
\begin{equation*}
S_{h,t}(Q)=tA(h)N(Q) +O_{h,\varepsilon} (Q^{3/2+\varepsilon}) .
\end{equation*}
\end{thm}

For $0<\alpha <2$, we define
\begin{equation*}
B_\alpha =\iint\limits_\TT \left[ \frac{1+s}{t} \right]^{\alpha}
ds\, dt =\sum_{k=1}^\infty
k^\alpha \ar (\TT_k)\ll \sum_{k=1}^\infty k^{\alpha-3} <\infty .
\end{equation*}
Employing the results from \cite{BCZ}, we give the
following generalization of \eqref{I7}:

\medskip

\begin{thm}\label{ThmB}
{\em (i)} For every $\alpha \in (0,2)$, we have
\begin{equation*}
\left| \sum\limits_{i} \nu_Q (\gamma_i)^\alpha -
2N(Q) B_\alpha \right| \ll_\alpha E_\alpha (Q)=\begin{cases}
Q\log Q & \mbox{if $\alpha <1$,} \\
Q\log^2 Q &\mbox{if $\alpha=1$,} \\
Q^\alpha \log Q & \mbox{if $1<\alpha <2$.}
\end{cases}
\end{equation*}

{\em (ii)} For every $\alpha \in ( 0,3/2 )$ and
$t\in [0,1]$, we have for all $\varepsilon >0$,
\begin{equation*}
\left| \sum\limits_{\gamma_i \leq t} \nu(\gamma_i)^\alpha -
2tN(Q) B_\alpha \right| \ll_{\alpha,\varepsilon} F_\alpha (Q)
=\begin{cases}
Q^{3/2+\varepsilon} & \mbox{if $\alpha \leq 1$,} \\
Q^{\alpha+1/2+\varepsilon} &\mbox{if $\alpha >1$.}
\end{cases}
\end{equation*}
\end{thm}

Note that
\begin{equation*}
B_1 =\sum\limits_{k=1}^\infty k\ar (\TT_k)=\frac{1}{6}+
\sum\limits_{k=2}^\infty \frac{4k}{k(k+1)(k+2)} =\frac{3}{2} \, ,
\end{equation*}
which is consistent with \eqref{I5} and \eqref{I4}.

Finally, we show that \eqref{I6} with error $O( Q\log Q/k)$
can be derived from Lemma 2 in \cite{BCZ}. In our framework the
geometrical significance of the constants
$l_k$ and $u_k$ is apparent, as
\begin{equation*}
l_k=2\ar (\TT_k^{\prime \prime})\qquad \mbox{\rm and} \qquad
u_k =2\ar (\TT_k^\prime).
\end{equation*}

Furthermore, if we set for $t\in (0,1]$,
\begin{equation*}
L(Q,k,t)=\# \bigg\{ \gamma_i=\frac{a_i}{q_i} \leq t\, ;\,
\nu_Q(\gamma_i)=k=\bigg[ \frac{2Q+1}{q_i}\bigg]-1 \bigg\}
\end{equation*}
and
\begin{equation*}
U(Q,k,t)=\# \bigg\{ \gamma_i=\frac{a_i}{q_i} \leq t\, ;\,
\nu_Q(\gamma_i) =k=\bigg[ \frac{2Q+1}{q_i}\bigg] \bigg\},
\end{equation*}
then we deduce as a result of Lemma 10 in \cite{BCZ} the
following generalization of \eqref{I6}:

\medskip

\begin{thm}\label{ThmC}
For every $t\in (0,1]$ and every $\varepsilon >0$, we have
\begin{equation*}
L(Q,k,t)=tl_k N(Q)+O_\varepsilon \bigg(
\frac{Q^{3/2+\varepsilon}}{k} \bigg)
\end{equation*}
and
\begin{equation*}
U(Q,k,t)=tu_k N(Q)+O_\varepsilon \bigg(
\frac{Q^{3/2+\varepsilon}}{k} \bigg) .
\end{equation*}
\end{thm}

\bigskip

{\bf Acknowledgements} We are grateful to Richard Hall and
Peter Shiu for correspondence on \cite{HS} and to the anonymous
referee for suggestions leading to the improvement of this paper.
\bigskip

\bigskip

\section{Proof of the main results}

We denote throughout
\begin{equation*}
\Zf=\{ (a,b)\in \Z^2 \, ;\, \gcd (a,b)=1\} .
\end{equation*}

By \eqref{I2}, we get
\begin{equation*}
\nu_Q (\gamma_{i+h} )=\bigg[ \frac{Q+q_{i+h-1}}{q_{i+h}} \bigg] =
\kappa_{h+1} \bigg( \frac{q_{i-1}}{Q},\frac{q_i}{Q} \bigg) ,
\end{equation*}
leading to
\begin{equation*}
S_h (Q)=\sum\limits_{\gamma_i \in \FF_Q} \kappa_1
\bigg( \frac{q_{i-1}}{Q},\frac{q_i}{Q} \bigg)
\kappa_{h+1} \bigg( \frac{q_{i-1}}{Q},\frac{q_i}{Q} \bigg) .
\end{equation*}

Since the pairs of denominators of consecutive elements in
$\FF_Q$ coincide with the elements of the set
\begin{equation*}
Q\TT \cap \Zf =\{ (a,b)\in \Zf \, ;\, a+b>Q\geq a,b\geq 1\} ,
\end{equation*}
we may write
\begin{equation*}
S_h(Q)=\sum\limits_{k=1}^\infty \sum\limits_{(a,b) \in
Q\TT \cap \Zf} \kappa_1 \bigg( \frac{a}{Q},\frac{b}{Q} \bigg)
\kappa_{h+1} \bigg( \frac{a}{Q},\frac{b}{Q} \bigg) .
\end{equation*}

Taking also stock on \eqref{I3}, this further yields
\begin{equation}\label{2.1}
\begin{split}
S_h (Q) & =\sum\limits_{k=1}^\infty
\sum\limits_{(a,b)\in Q\TT_k \cap \Zf}
\kappa_{h+1} \bigg( \frac{a}{Q},\frac{b}{Q} \bigg) \\
& =\sum\limits_{k=1}^\infty k \sum\limits_{l=1}^\infty
l\, \# \bigg\{ (a,b)\in Q\TT_k \cap \Zf \, ;\,
\kappa_1 T^h \bigg( \frac{a}{Q},\frac{b}{Q} \bigg)=l\bigg\} \\
& =\sum\limits_{k,l=1}^\infty kl\, \# \bigg\{
(a,b)\in Q\TT_k \cap \Zf \, ;\, T^h \bigg( \frac{a}{Q},\frac{b}{Q}
\bigg) \in \TT_l \bigg\} \\
& =\sum\limits_{k,l=1}^\infty kl\, \# \big( Q(\TT_k \cap
T^{-h} \TT_l )\cap \Zf \big) .
\end{split}
\end{equation}

If we set
\begin{equation*}
\TT_k^* =\bigcup\limits_{n=k}^\infty \TT_n ,
\end{equation*}
then
\begin{equation}\label{2.2}
\TT_k^* \cap T^{-h} \TT_l^* =\bigg( \bigcup\limits_{m=k}^\infty
\TT_m \bigg) \cap T^{-h} \bigg( \bigcup\limits_{n=l}^\infty \TT_n
\bigg) =\bigcup\limits_{m=k}^\infty \bigcup\limits_{n=l}^\infty
( \TT_m \cap T^{-h} \TT_n ) .
\end{equation}

We wish to estimate
\begin{equation*}
A_{k,l}(Q)=\# \big( Q(\TT_k^* \cap T^{-h} \TT_l )\cap \Zf \big) .
\end{equation*}

Since the sets from the right-hand side of \eqref{2.2} are
mutually disjoint, we now have
\begin{equation*}
\begin{split}
\sum\limits_{k,l=1}^\infty A_{k,l}(Q) &  =\sum\limits_{k,l=1}^\infty
\sum\limits_{m=k}^\infty \sum\limits_{n=l}^\infty \# \big(
Q(\TT_m \cap T^{-h} \TT_n )\cap \Zf \big)  \\
& =\sum\limits_{m,n=1}^\infty \#
\big( Q(\TT_m \cap T^{-h} \TT_n )\cap \Zf \big)
\sum\limits_{k=1}^m \sum\limits_{l=1}^n 1 \\ &
=\sum\limits_{m,n=1}^\infty mn\, \# \big( Q(\TT_m \cap
T^{-h} \TT_n )\cap \Zf \big) ,
\end{split}
\end{equation*}
so that, by \eqref{2.1},
\begin{equation}\label{2.3}
S_h(Q) =\sum\limits_{k,l=1}^\infty A_{k,l} (Q) .
\end{equation}

The following lemma is proved in a similar way as
Lemma 2 in \cite{BCZ}.

\medskip

\begin{lem}\label{L1}
Let $\Omega \subset [0,R_1] \times [0,R_2 ]$ be a region in
$\R^2$ with rectifiable boundary $\partial \Omega$ and
let $g:\Omega \rightarrow \R$ be a $C^1$ function on $\Omega$.
Suppose that $R\geq \min (R_1,R_2)$. Then we have
\begin{equation*}
\begin{split}
\sum\limits_{(a,b)\in \Omega \cap \Zf} \hspace{-8pt}
g(a,b) & =\frac{6}{\pi^2} \iint\limits_\Omega g(x,y)\, dxdy
+O\big( \| Dg\|_\infty \ar (\Omega) \log R\big) \\ & +
O\bigg( \| g\|_\infty \Big( R+\frac{\ar (\Omega)}{R}+
\lh (\partial \Omega)\log R \Big) \bigg) ,
\end{split}
\end{equation*}
where
\begin{equation*}
\| Dg\|_\infty =\sup\limits_{(x,y)\in \Omega}
\bigg| \frac{\partial g}{\partial x} (x,y) \bigg|
+ \bigg| \frac{\partial g}{\partial y} (x,y) \bigg| .
\end{equation*}
\end{lem}

\medskip

\begin{cor}\label{C2}
Let $\Omega \subseteq [0,R_1]\times [0,R_2]$ be a bounded region
with rectifiable boundary $\partial \Omega$ and let
$R\geq\min (R_1,R_2 )$. Then we have
\begin{equation*}
\# (\Omega \cap \Zf)=\frac{6\ar(\Omega)}{\pi^2}+
O\left( R+\lh(\partial \Omega)\log R+\frac{\ar(\Omega)}{R}\right).
\end{equation*}
\end{cor}

\smallskip

\begin{remark}\label{R3}
Let $\left( \begin{smallmatrix} a & b \\ c & d \end{smallmatrix}
\right) \in SL_2 (\Z)$, let $\Phi$ be the linear transformation
on $\R^2$ defined by
$\Phi(x,y)=(ax+by,cx+dy)$, and let $\Omega$ be a bounded region
in $\R^2$. Then $\# (\Omega \cap \Zf)=\# (\Phi \Omega \cap \Zf)$,
which implies in turn that
\begin{equation*}
\# (Q\Omega \cap \Zf)=
\# \big( Q(T\Omega )\cap \Zf \big) =\dots =
\# \big( Q(T^h \Omega )\cap \Zf \big)
\end{equation*}
whenever $\Omega \subseteq \TT$. In particular, this gives
\begin{equation}\label{2.4}
\# \big( Q(\TT_k^* \cap T^{-h} \TT_l^* )\cap \Zf \big) =
\# \big( Q(T^h \TT_k^* \cap \TT_l^* )\cap \Zf \big).
\end{equation}
\end{remark}

\medskip

\begin{lem}\label{L4}
{\em \cite[Lemma\,5]{BCZ}}
Let $r\geq 1$ and suppose that $i\neq j$ and that
$\max \big( L_i(x,y),L_j(x,y)\big) \leq 2^{-r-1}$ for some point
$(x,y)\in \TT$. Then $\vert j-i\vert >r+1$.
\end{lem}

\smallskip

\begin{cor}\label{C5}
Assume that $\min (k,l)>c_h =2^{h+1}$. Then we have
\begin{equation*}
\TT_k^* \cap T^{-h} \TT_l^* =\emptyset .
\end{equation*}
\end{cor}

{\sl Proof.} Suppose that there exists $(x,y)\in \TT_k^*
\cap T^{-h} \TT_l^*$. Then $L_1 (x,y)=y\leq 2/k <2^{-h}$ and
$L_1 \big( T^h (x,y)\big) =L_{h+1} (x,y)\leq 2/l <
2^{-h}$. We now infer from Lemma \ref{L4} that $\vert h+1-1\vert
>h$, which is a contradiction. \qed

\smallskip

\begin{remark}\label{R6}
A careful look at the proof of Lemma \ref{L4} shows that the constant
$2^{r+1}$ can be lowered to $4r+2$ (see Lemma 3.4 and Remark 3.5 in
\cite{BCZ1}). Note also that $T^{\pm 1} \TT_m^* \subset \TT_1$ for
all $m\geq 5$. By Lemma 3.4 in \cite{BCZ1} it also follows that
for all $r\geq 2$ and all $m\geq c_r=4r+2$, we have
$\cup_{i=2}^r T^{\pm i} \TT_m^* \subset \TT_2$. In summary, we have
\begin{equation*}
T^h \TT_m^* \cap \TT_n^* =\emptyset
\end{equation*}
whenever $m\geq c_h$ and $n\geq 3$.
\end{remark}

Since $Q\TT_k^* \cap \Zf =\emptyset$ whenever $2Q/k<1$,
Remark \ref{R3} shows that $A_{k,l}(Q)=0$ unless
$\max (k,l)\leq 2Q$. By Corollary \ref{C5}, $A_{k,l}(Q)=0$
unless $\min (k,l)\leq c_h$. Thus \eqref{2.3} provides
\begin{equation}\label{2.5}
S_h(Q)=\sum\limits_{k,l=1}^{2Q} A_{k,l}(Q)=C_h(Q)+D_h(Q),
\end{equation}
where
\begin{equation*}
C_h(Q)=\sum\limits_{k=c_h+1}^{2Q} \sum\limits_{l=1}^{c_h}
A_{k,l} (Q)\qquad \mbox{\rm and} \qquad
 D_h(Q)=\sum\limits_{k=1}^{c_h} \sum\limits_{l=1}^{2Q} A_{k,l} (Q).
\end{equation*}

Suppose first that $k>c_h$. Applying Corollary \ref{C2} to
$\Omega =Q(\TT_k^* \cap T^{-h} \TT_l^*)$ with
$\ar (\Omega )\ll Q^2/k^2$, and with
\begin{equation*}
\begin{split}
\lh (\partial \Omega ) & =Q\lh \bigg( \bigcup\limits_{r=l}^\infty
T^{-h} \TT_r \cap \TT_k^* \bigg) =Q\lh \bigg(
\bigcup\limits_{r=l}^{c_h} T^{-h} \TT_r \cap \TT_k^* \bigg) \\
& \leq Q\sum\limits_{r=l}^{c_h} \lh (T^{-h} \TT_r \cap \TT_k^* )
\ll_h \frac{Q}{k} \, ,
\end{split}
\end{equation*}
we infer that
\begin{equation}\label{2.6}
A_{k,l}(Q)=\frac{6Q^2}{\pi^2} \ar(\TT_k^* \cap T^{-h} \TT_l^*)+
O_h \left( \frac{Q\log Q}{k} \right) .
\end{equation}
\medskip

By \eqref{2.6} and $\ar (\TT_k^* )=O(k^{-2})$ we gather
\begin{equation}\label{2.7}
\begin{split}
C_h(Q) & =\frac{6Q^2}{\pi^2} \sum\limits_{k=c_h+1}^{2Q}
\sum\limits_{l=1}^{c_h} \ar (\TT_k^* \cap T^{-h} \TT_l^*)
+O_h (Q\log^2 Q) \\
& =\frac{6Q^2}{\pi^2} \sum\limits_{k=c_h+1}^\infty
\sum\limits_{l=1}^\infty
\ar (\TT_k^* \cap T^{-h} \TT_l^*)+O_h (Q\log^2 Q).
\end{split}
\end{equation}

When $k\leq c_h$, we employ Remark \ref{R3} and equality
\eqref{2.4} to get
\begin{equation*}
A_{k,l}(Q)=\# \big( Q (\TT_l^* \cap T^h \TT_k^* )\big) .
\end{equation*}
Applying now Corollary \ref{C2} to $\Omega =Q(\TT_l^* \cap
T^h \TT_k^*)$ with $\lh (\partial \Omega)
\ll_h Q/l$, and using the fact that $T$ is
area-preserving, we obtain
\begin{equation}\label{2.8}
\begin{split}
A_{k,l}(Q) & =\frac{6}{\pi^2} \, \ar (\TT_l^* \cap
T^h \TT_k^* )+O_h \bigg( \frac{Q\log Q}{l} \bigg) \\ &
=\frac{6}{\pi^2} \, \ar (\TT_k^* \cap T^{-h} \TT_l^*)
+O_h \bigg( \frac{Q\log Q}{l} \bigg) .
\end{split}
\end{equation}
We may now employ \eqref{2.7} and $\ar (\TT_l^*)=O(l^{-2})$, to get
\begin{equation}\label{2.9}
D_h(Q)=\frac{6}{\pi^2} \sum\limits_{k=1}^{c_h}
\sum\limits_{l=1}^\infty \ar (\TT_k^* \cap T^{-h} \TT_l^* )+
O_h (Q\log^2 Q).
\end{equation}

Inserting \eqref{2.9} and \eqref{2.7} into \eqref{2.5}, we collect
\begin{equation}\label{2.10}
S_h(Q)=\frac{6}{\pi^2} \sum\limits_{k,l=1}^\infty \ar
(\TT_k^* \cap T^{-h} \TT_l^*) +O_h (Q\log^2 Q).
\end{equation}

We conclude the proof of Theorem \ref{ThmA} (i) by noting as in
\eqref{2.3} that the sum in the right-hand side of \eqref{2.10}
is equal to
\begin{equation}\label{2.11}
\sum\limits_{m,n=1}^\infty \ar (T^h \TT_m^* \cap \TT_n^* )=
\sum\limits_{m,n=1}^\infty mn\ar (T^h \TT_m \cap \TT_n)=
\frac{A(h)}{2} \, .
\end{equation}

We set
\begin{equation*}
f=f_0+R_0 \qquad \mbox{\rm and} \qquad
g=f\circ T^{-h} =f_0 \circ T^{-h}+R_0 \circ T^{-h} ,
\end{equation*}
with
\begin{equation*}
f_0 =\sum\limits_{m=1}^{c_h-1} e_{\TT_m^*},\qquad
R_0 =\sum\limits_{m=c_h}^\infty e_{\TT_m^*} ,
\end{equation*}
where $e_S$ denotes the characteristic function of the set $S$.
By Remark \ref{R6} we gather that the product of $R_0$ and $g$ equals
\begin{equation*}
\sum\limits_{m=c_h}^\infty \sum\limits_{n=1}^\infty
e_{T^h \TT_m^* \cap \TT_n^*} =\sum\limits_{m=c_h}^\infty
\sum\limits_{n=1}^2 e_{T^h \TT_m^* \cap \TT_n^*} .
\end{equation*}
Hence,
\begin{equation*}
\iint\limits_\TT R_0 (s,t) g(s,t)\, ds\, dt \leq
2\sum\limits_{m=c_h}^\infty \ar (\TT_m^*) \ll \sum\limits_{m=1}^\infty
\frac{1}{m^2} \ll 1.
\end{equation*}
In a similar way, we get
\begin{equation*}
\iint\limits_\TT f(s,t) (R_0 \circ T^{-h})(s,t) \, ds\, dt \ll 1.
\end{equation*}

Using \eqref{2.11} we may now write
\begin{equation}\label{2.12}
\begin{split}
A(h) & =2\iint\limits_\TT f(s,t)(f\circ T^{-h})(s,t)\, ds\, dt
\\ & \ll 1+2\iint\limits_\TT f_0(s,t)g_0(s,t)\, ds \, dt.
\end{split}
\end{equation}
The function $f_0$ belongs to $L^2 (\TT)$ as a result of
\begin{equation*}
\begin{split}
f_0 & =\sum\limits_{m=1}^{c_h-1} \sum\limits_{k=m}^\infty e_{\TT_k} =
\sum\limits_{k=1}^\infty \sum\limits_{m=1}^{\min (k,c_h-1)}
e_{\TT_k} =\sum\limits_{k=1}^\infty \min (k,c_h-1) e_{\TT_k} \\
& =\sum\limits_{k=1}^{c_h-1} ke_{\TT_k} +(c_h-1)
\sum\limits_{k=c_h}^\infty e_{\TT_k} =
\sum\limits_{k=1}^{c_h-1} ke_{\TT_k} +(c_h-1) e_{\TT_{c_h}^*} ,
\end{split}
\end{equation*}
which gives in turn
\begin{equation*}
\| f_0 \|_{L^2(\TT)} =\sum\limits_{k=1}^{c_h-1} k^2 \ar (\TT_k) +
(c_h-1)^2 \ar (\TT_{c_h}^*) \ll \sum\limits_{k=1}^{c_h}
\frac{1}{k}+1 \ll 1+\log h.
\end{equation*}
Since $T$ is area-preserving, we also have
\begin{equation*}
\| g_0 \|_{L^2 (\TT)} =\| f_0 \circ T^{-h} \|_{L^2 (\TT)} =
\| f_0 \|_{L^2(\TT)} ,
\end{equation*}
and the Cauchy-Schwarz inequality gives in conjunction
with \eqref{2.12}
\begin{equation*}
A(h) \ll 1+2\| f_0 \|_{L^2 (\TT)} \| g_0 \|_{L^2 (\TT)}
\ll 1+\log h.
\end{equation*}

To prove (ii), we proceed as in Section 8 in \cite{BCZ} and note
that the relation $a^\prime q-aq^\prime =1$ between two
consecutive Farey fractions of order $Q$ shows that
$a=q-\bar{q^\prime}$, where $\bar{q^\prime}
\in \{ 1,\dots,q\}$ denotes the multiplicative inverse of
$q^\prime \!\!\pmod{q}$. Thus the condition $a/q \in I$
with $I=(\alpha,\beta ]\subseteq (0,1]$ interval, is equivalent to
$\bar{q^\prime} \in I_q =[ (1-\beta)q,(1-\alpha)q )$.
As a result, we may write
\begin{equation*}
S_{h,t}(Q)=\sum_{\substack{(a,b)\in Q\TT \\ \bar{b}\in I_a}}
\kappa_1 \bigg( \frac{a}{Q},\frac{b}{Q} \bigg)
\kappa_{h+1} \bigg( \frac{a}{Q},
\frac{b}{Q} \bigg) =\sum\limits_{k,l=1}^\infty A_{k,l,t} (Q),
\end{equation*}
with
\begin{equation*}
A_{k,l,t} (Q)=\# \big( Q(\TT_k^* \cap T^{-h} \TT_l^* ) \cap
\{ (a,b)\, ;\, \bar{b} \in I_a \} \big),
\end{equation*}
where $I=(0,t]$. Applying Lemma 10 in \cite{BCZ} to $f(a,b)=1$,
$\Omega =Q(\TT_k^* \cap T^{-h} \TT_l^* )$, $A=O(Q)$,
$R_1=R_2 =O( Q/k )$, and employing \eqref{2.6}
and \eqref{2.8}, we infer that
\begin{equation*}
\begin{split}
A_{k,l,t}(Q) & =tA_{k,l}(Q)+O\left( \frac{Q}{k^2}+
\frac{Q^{3/2+\varepsilon}}{k} \right) \\ & =
\frac{6tQ^2}{\pi^2} \, \ar (\TT_k^* \cap T^{-h} \TT_l^*)
+O\left( \frac{Q^{3/2+\varepsilon}}{k} \right),
\end{split}
\end{equation*}
leading to
\begin{equation*}
S_{h,t}(Q)=\frac{3tA(h)Q^2}{\pi^2}+O(Q^{3/2+\varepsilon}\log Q).
\end{equation*}
\qed

\medskip

{\sl Proof of Theorem {\em \ref{ThmB}}}.
This is similar to the proof of Theorem \ref{ThmA}. We get as in
\eqref{2.3}
\begin{equation}\label{2.13}
\begin{split}
T_\alpha (Q) & :=\sum\limits_i \nu_Q (\gamma_i)^\alpha =
\sum\limits_{k=1}^\infty k^\alpha \, \# (Q\TT_k \cap \Zf ) \\ &
=\sum\limits_{k=1}^\infty \# (Q\TT_k \cap \Zf ) \sum\limits_{m=1}^k
\big( m^\alpha -(m-1)^\alpha \big) \\
& =\sum\limits_{m=1}^\infty \big( m^\alpha -(m-1)^\alpha \big)
\sum\limits_{k=m}^\infty \# (Q\TT_k \cap \Zf ) \\ &
=\sum\limits_{m=1}^{2Q} \big( m^\alpha -(m-1)^\alpha \big)
\, \# (Q\TT_m^* \cap \Zf ).
\end{split}
\end{equation}

However, Corollary \ref{C2} gives
\begin{equation*}
\# (Q\TT_m^* \cap \Zf )=\frac{6Q^2}{\pi^2} \,
\ar (\TT_m^*)+O\left( \frac{Q\log Q}{m} \right) ,
\end{equation*}
which we insert into \eqref{2.13} to get
\begin{equation*}
T_\alpha (Q)=\frac{6Q^2}{\pi^2} \sum\limits_{m=1}^{2Q}
\big( m^\alpha -(m-1)^\alpha \big) \ar (\TT_m^* )+
O_\alpha \big( E_\alpha (Q)\big),
\end{equation*}
with $E_\alpha (Q)$ as in Theorem \ref{ThmB}.

However, $Q^2 \sum_{m=Q}^\infty m^{\alpha -1} \ar (\TT_m^*) \ll Q^2
\sum_{m=Q}^\infty m^{\alpha-3} \ll Q^\alpha$. Hence,
\begin{equation}\label{2.14}
T_\alpha (Q)=\frac{6Q^2}{\pi^2} \sum\limits_{m=1}^\infty
\big( m^\alpha -(m-1)^\alpha \big) \ar (\TT_m^*) +O_\alpha
\big( E_\alpha (Q)\big) .
\end{equation}

Finally, Theorem \ref{ThmB} (i) follows from \eqref{2.14}, and from
\begin{equation}\label{2.15}
\begin{split}
\sum\limits_{m=1}^\infty & \big( m^\alpha -(m-1)^\alpha \big)
\ar (\TT_m^* ) \\ & =
\sum\limits_{m=1}^\infty \big( m^\alpha -(m-1)^\alpha \big)
\sum\limits_{k=m}^\infty \ar (\TT_k) \\ &  =
\sum\limits_{k=1}^\infty \ar (\TT_k) \sum\limits_{m=1}^k
\big( m^\alpha -(m-1)^\alpha \big)
= \sum\limits_{k=1}^\infty k^\alpha \ar (\TT_k) \\ & =
\iint\limits_\TT \Big[ \frac{1+s}{t} \Big]^\alpha \, ds\, dt.
\end{split}
\end{equation}

To prove Theorem \ref{ThmB} (ii), we first write as in \eqref{2.13}
\begin{equation}\label{2.16}
T_{\alpha,t} (Q) :=\sum\limits_{\gamma_i \leq t}
\nu_Q (\gamma_i)^\alpha =\sum\limits_{m=1}^{2Q}
\big( m^\alpha -(m-1)^\alpha \big)
\sum_{\substack{(a,b)\in Q\TT_m^* \\ \bar{b} \in I_a}} 1,
\end{equation}
where $I=(0,t]$ and $I_a =[ (1-t)a,a)$. Applying
Lemma 10 in \cite{BCZ} to $f(a,b)=1$ and $\Omega =Q\TT_m^*$, we find
\begin{equation*}
\sum_{\substack{(a,b)\in Q\TT_m^* \\ \bar{b} \in I_a}} 1 =
\frac{6tQ^2}{\pi^2} \, \ar (\TT_m^*) +
O\left( \frac{Q^{3/2+\varepsilon}}{m} \right),
\end{equation*}
which gives in conjunction with \eqref{2.16} and \eqref{2.15}
\begin{equation*}
\begin{split}
T_{\alpha,t} (Q) & =\frac{6tQ^2}{\pi^2} \sum\limits_{m=1}^{2Q}
\big( m^\alpha -(m-1)^\alpha \big) \ar (\TT_m^*) +O_\alpha
\big( F_\alpha (Q)\big) \\ & =\frac{6tQ^2}{\pi^2}
\iint\limits_\TT \Big[ \frac{1+s}{t} \Big]^\alpha \, ds\, dt +
O_\alpha \big( F_\alpha (Q)\big),
\end{split}
\end{equation*}
where $F_\alpha (Q)$ is as in the statement of Theorem \ref{ThmB}.

The equalities
\begin{equation*}
\nu_Q (\gamma_i)=\bigg[ \frac{Q+q_{i-1}}{q_i} \bigg] =k=
\bigg[ \frac{2Q+1}{q_i} \bigg] -1
\end{equation*}
read as
\begin{equation*}
(q_{i-1},q_i)\in Q\TT_k^L \cap \Zf,
\end{equation*}
where
\begin{equation*}
\TT_k^L =\left\{ (s,t)\in \TT_k^L \, ;\, \frac{2+1/Q}{k+2}
<t\leq \frac{2+1/Q}{k+1} \right\} .
\end{equation*}
Hence $L(Q,k)=\# (Q\TT_k^L \cap \Zf)$ and by means of Corollary
\ref{C2}, we find
\begin{equation*}
\begin{split}
L(Q,k) & =\frac{6Q^2}{\pi^2} \, \ar (\TT_k^L)+
O\bigg( \frac{Q}{k}+\frac{Q\log Q}{k}+\frac{Q}{k^2} \bigg) \\ & =
\frac{6Q^2}{\pi^2} \, \ar (\TT_k^{\prime \prime})+
O\bigg( \frac{Q\log Q}{k} \bigg) .
\end{split}
\end{equation*}

The estimate on $U(Q,k)$ in \eqref{I6} follows in a similar way.

Theorem \ref{ThmC} is derived from Lemma 10 in \cite{BCZ} in
a similar way as above.

\bigskip

\section{Some numerical computations}

The transformation $T$ maps each region $\TT_k$ onto its symmetric
with respect to the first bisector. That is,
\begin{equation}\label{3.1}
T\TT_k =S\TT_k \qquad \mbox{\rm and} \qquad T\TT_k^* =S\TT_k^* ,
\end{equation}
where $S$ acts on $\TT$ as
\begin{equation*}
S(x,y)=(y,x) .
\end{equation*}
Moreover, we notice that the inverse $T^{-1}$ of $T$ can be
expressed as
\begin{equation}\label{3.2}
T^{-1}=STS.
\end{equation}
We now infer from \eqref{3.2}, \eqref{3.1}, and the fact that $T$
is area-preserving that
\begin{equation}\label{3.3}
\begin{split}
\ar (T^h \TT_m^* \cap \TT_n^* ) &
=\ar (\TT_m^* \cap T^{-h} \TT_n^* )
=\ar (\TT_m^* \cap ST^h S\TT_n^* ) \\ &
=\ar (S\TT_m^* \cap T^h S\TT_n^* ) =
\ar (T\TT_m^* \cap T^h T\TT_n^* ) \\ &
=\ar (T^h \TT_n^* \cap \TT_m^* ).
\end{split}
\end{equation}
Thus, to evaluate $A(h)$ via \eqref{2.11},
it suffices to consider $m\leq n$ only.
We find the following table for the value of
$\ar (T\TT_m^* \cap \TT_n^*)$
{\small $$\begin{tabular}{|c||c|c|c|c|c|} \hline
 & $n=1$ & $n=2$ & $n=3$ & $n=4$ & $n\geq 5$ \\ \hline \hline
$m=1$ & $\frac{1}{2}$ & $\frac{1}{3}$ & $\frac{1}{6}$ &
$\frac{1}{10}$ & $\frac{2}{n(n+1)}$ \\ \hline
$m=2$ & $\frac{1}{3}$ & $\frac{1}{6}$ & $\frac{1}{30}$ &
$\frac{1}{210}$ & $0$ \\ \hline
$m=3$ & $\frac{1}{6}$ & $\frac{1}{30}$ & $0$ & $0$ & $0$ \\ \hline
$m=4$ & $\frac{1}{10}$ & $\frac{1}{210}$ & $0$ & $0$ & $0$ \\ \hline
$m\geq 5$ & $\frac{2}{m(m+1)}$ & $0$ & $0$ & $0$ & $0$ \\
\hline
\end{tabular}
$$}
which gives in turn
\begin{equation*}
A(1)=2\cdot \frac{96}{35} =\frac{192}{35} \approx 5.4857.
\end{equation*}

We also find the following table for the value of
$\ar (T^2 \TT_m^* \cap \TT_n^* )$, which is symmetric as a
result of \eqref{3.3},
{\small $$\begin{tabular}{|c||c|c|c|c|c|c|c|c|c|} \hline
 & $n=1$ & $n=2$ & $n=3$ & $n=4$ & $n=5$ &
$n=6$ & $n=7$ & $n=8$ & $n\geq 9$ \\ \hline \hline
$m=1$ & $\frac{1}{2}$ & $\frac{1}{3}$ & $\frac{1}{6}$ &
$\frac{1}{10}$ & $\frac{1}{15}$ & $\frac{1}{21}$ & $\frac{1}{28}$ &
$\frac{1}{36}$ & $\frac{2}{n(n+1)}$ \\ \hline
$m=2$ & $\frac{1}{3}$ & $\frac{23}{84}$ & $\frac{31}{210}$ &
$\frac{2}{21}$ & $\frac{1}{15}$ & $\frac{1}{21}$ & $\frac{1}{28}$
& $\frac{1}{36}$ & $\frac{2}{n(n+1)}$  \\ \hline
$m=3$ & $\frac{1}{6}$ & $\frac{31}{210}$ & $\frac{1}{10}$ &
$\frac{13}{210}$ & $\frac{1}{30}$ & $\frac{1}{70}$ & $\frac{1}{220}$ &
$\frac{1}{1170}$ & $0$ \\ \hline
$m=4$ & $\frac{1}{10}$ & $\frac{2}{21}$ & $\frac{13}{210}$ &
$\frac{1}{42}$ & $\frac{1}{231}$ & $0$ & $0$ & $0$ & $0$ \\ \hline
$m=5$ & $\frac{1}{15}$ & $\frac{1}{15}$ & $\frac{1}{30}$ &
$\frac{1}{231}$ & $0$ & $0$ & $0$ & $0$ & $0$ \\ \hline
$m=6$ & $\frac{1}{21}$ & $\frac{1}{21}$ & $\frac{1}{70}$ &
$0$ & $0$ & $0$ & $0$ & $0$ & $0$ \\ \hline
$m=7$ & $\frac{1}{28}$ & $\frac{1}{28}$ & $\frac{1}{220}$ &
$0$ & $0$ & $0$ & $0$ & $0$ & $0$ \\ \hline
$m=8$ & $\frac{1}{36}$ & $\frac{1}{36}$ & $\frac{1}{1170}$ &
$0$ & $0$ & $0$ & $0$ & $0$ & $0$ \\ \hline
$m\geq 9$ & $\frac{2}{m(m+1)}$ & $\frac{2}{m(m+1)}$ & $0$ & $0$ & $0$ &
$0$ & $0$ & $0$ & $0$ \\ \hline
\end{tabular}$$}
and collect
\begin{equation*}
A(2)=2\cdot \frac{796727}{180180}=\frac{796727}{90090}
\approx 8.8437.
\end{equation*}

Using the second part of Remark \ref{R6} and the fact that
$T$ is area-preserving, we may write for $h\geq 2$,
\begin{equation*}
\begin{split}
\frac{A(h)}{2} & =\sum\limits_{m,n=1}^{c_h-1} mn\ar (T^h \TT_m \cap
\TT_n) +4\sum\limits_{n=c_h}^\infty n\ar (\TT_n )  \\
& =\sum\limits_{m,n=1}^{c_h-1} mn\ar (T^h \TT_m \cap
\TT_n)+\frac{2}{(c_h+1)(c_h+2)}\, .
\end{split}
\end{equation*}
However, each region $T^r \TT_m$ is a finite union of triangles
with rational numbers as vertex coordinates. Thus $A(h)$ is
a rational number for any $h\in \N^*$.

\bigskip

\bigskip

\bibliographystyle{amsplain}

\end{document}